\documentstyle{article}

\begin{document}

\baselineskip = 15 pt

\title{\bf $\hbar$-(Yangian) Deformation of Miura Map \\
and Virasoro Algebra}

\vspace{1cm}
\author{Xiang-Mao Ding\\
Institute of Theoretical Physics,
Academy of China, Beijing 100080, China\\
 \vspace{2mm}\\
Bo-Yu Hou \hspace{1cm}  Liu Zhao\\
Institute of Modern Physics, Northwest University, Xian 710069,
China}

\date{December 1996}
\maketitle

\begin{abstract}
An $\hbar$-deformed Virasoro Poisson algebra is obtained using
the Wakimoto realization of the Sugawara operator for the Yangian double
$DY_\hbar(sl_2)_c$ at the critical level $c=-2$.
\end{abstract}

\section{Introduction}
In this work we construct the Yangian deformed Miura map and the
corresponding (deformed) Virasoro algebra. Since Drinfeld
proposed in his works \cite{D1,D2}, the Yangian algebra $Y_\hbar(g)$ is
known as a kind of Hopf algebra associated with rational solutions of
quantum Yang-Baxter equation characterized by an additive
parameter $\hbar$, and is a deformation of the universal enveloping algebra
of the loop algebra $g[u,u^{-1}]$ associated with the
simple Lie algebra $g$. Since the rational solution of the
Yang-Baxter equation is characterized
by an additive parameter $\hbar$, we call the corresponding deformed
algebras $\hbar$-deformations. Many faces of the Yangian algebra have been
thoroughly studied in the literatures with emphasis in both algebraic
aspects and physical applications. Algebraically, the development of
studies of Yangian algebra is almost parallel to the development of
studies of the so-called quantum algebra associated with
trigonometric solutions of Yang-Baxter solution and characterized by a
multiplicative parameter $q$ (and hence the
name $q$-deformed algebras \cite{D3}),
and such parallelism became even more perfect after the work of Khoroshkin
{\it et al} who successfully constructed the quantum double for
Yangian algebra (the Yangian double \cite{KT})
and made a central extension for Yangian double \cite{K,KL,KLP,IK}.
Physically the Yangian algebra (actually its quantum double $DY_\hbar(g)$)
is found to be the dynamical symmetry algebra of many massive quantum
integrable field theories, among which we would like to mention the
Thirring model, principal chiral model and nonlinear sigma model
\cite{S2,BL,LS}. Generally speaking
the field theoretical models yielding Yangian symmetries often correspond
to certain scaling limits of some lattice statistical models away
from criticality. In this sense the Yangian algebra is in a
position somewhere in between the usual Kac-Moody Lie algebra and the
quantum algebra.

In the study of $q$-deformed Lie algebras there were a long standing
problem which is not resolved until recently, i.e. the construction of
a $q$-deformation of Virasoro algebra. After several attempts by different
authors (see, e.g. \cite{eg}), in Ref.\cite{FR}, Frenkel and
Reshetikhin obtained a version of $q$-deformed Virasoro and together
W algebras as $q$-deformed Gelfand-Dickey Poisson algebras. Later on,
the quantum version of their algebras has also been obtained by several
groups, see Refs. \cite{SKAO,LUKU,AKOS,AKOS2,FF}. Vertex operators
connected to such algebras have also been studied in Refs.
\cite{LuPu,SKMAO,AJMP}. The central idea of Ref. \cite{FR} can be
briefly summarized as follows. In the undeformed case, the (quantum)
Virasoro algebra can be constructed from the Kac-Moody algebra, e.g.
$\widehat{sl}_2$, by means of Sugawara construction, and by an appropriate
renormalization, the Virasoro generating function becomes a center
in the formal completion of the universal enveloping algebra
of $\widehat{sl}_2$ at the critical level $k=-2$
(the dual Coexeter number with a minus sign).
The Poisson brackets for the Virasoro algebra can then
be obtained from the Wakimoto realization of the $\widehat{sl}_2$ Kac-Moody
current in the limit $k+2 \rightarrow 0$, and this construction has a natural
connection to the famous Miura transformation (actually the Virasoro
Poisson brackets has to be obtained via this transformation). In the
$q$-deformed case, Frenkel and Reshetikhin \cite{FR}
successfully made a parallel development. Using the Ding-Frenkel
equivalence \cite{DF} of Drinfeld
currents  \cite{Dr:new} and the Reshetikhin-Semenov-Tian-Shansky
realization \cite{RS} of
$q$-affine algebras they obtained the center
of the formal completion of the $q$-algebra $U_q(\widehat{sl}_2)_k$. Then
Using the $q$-deformed Wakimoto realization \cite{AOS} of the Drinfeld
currents they showed that that center is nothing but a $q$-deformation of the
Miura map. Finally they obtained the Poisson bracket algebra for $q$-deformed
Virasoro algebra using the $q$-Miura map at the critical level $k=-2$. The
$q$-deformed W algebras are also obtained in a similar spirits.

It is interesting to ask that whether the constructions worked in
undeformed and $q$-deformed cases also works in $\hbar$-deformed case.
The answer is true but the construction is in some sense not straightforward,
as will be shown in the main context of this paper.
The $\hbar$-deformation of Virasoro algebra is our central object,
and we feel that this algebra is important to be studied in detail
because it is the corner stone of several important algebraic objects:
(i) it is a deformation of the conventional Virasoro algebra; (ii) it is
the scaling limit of the $q$-deformed Virasoro algebra obtained in
Ref. \cite{FR}; (iii) its connection with Yangian algebra with center
is precisely the sort of connections between the $q$-deformed Virasoro
and $q$-affine algebras; (iv) it is the classical counterpart of the
quantum $\hbar$-deformed Virasoro algebra obtained from the quantum
$q$-Virasoro algebra by taking the scaling limit \cite{YWL}. Moreover, the
connections to the $\hbar$-deformed Miura map is also an important
problem because that, while extended to algebras of higher rank,
this may reveal a new kind of deformed Gelfand-Dickey equation and may
also play some role in an $\hbar$-deformed Drinfeld-Sokolov reduction
scheme \cite{DSR}.

The outline of this work is as follows. In Section 2 we shall collect
necessary backgrounds and formulas by a brief review of the
Yangian doubles $DY({gl}_2)_c$ and $DY_\hbar({sl}_2)_c$. Section 3 is
devoted to the construction of $\hbar$-deformed Sugawara operator.
Then, using the Wakimoto realization given in \cite{KNO}, we derive the
$\hbar$-deformed Miura map in Section 4. In Section 5 we present the
Poisson bracket for the $\hbar$-deformed Virasoro algebra and Section 6
is devoted to some discussions and out-lookings.

\section{The Yangian algebras $DY_\hbar({gl}_2)_c$ and
$DY_\hbar({sl}_2)_c$}

The Yangian algebra we shall make use of is actually the central extension
$DY_\hbar( {gl}_2)_c$ and $DY_\hbar({sl}_2)_c$ of the quantum
doubles of $Y_\hbar({gl}_2)$ and $Y_\hbar(sl_2)$ respectively with
central element $c$. These algebras can be realized in three equivalent
ways, namely using the Chevalley generators, Drinfeld currents and
Reshetkhin-Semenov-Tian-Shansky formalism. In this section we shall collect
necessary formulas by making a brief review of the algebra
$DY_\hbar({gl}_2)_c$ and treating the algebra
$DY_\hbar({sl}_2)_c$ as the subalgebra of
$DY_\hbar({sl}_2)_c$ modulo a Heisenberg subalgebra.

In terms of Drinfeld currents, the algebra $DY_\hbar({gl}_2)_c$
can be regarded as the formal completion of the algebra generated
by the currents $k^{\pm}_i(u)$ ($i=1,2$), $e^\pm(u)$,
$f^{\pm}(u)$ together with a derivative $d$ and a center element $c$
with the following generating relations \cite{IK},

\begin{eqnarray}
& & [d,e(u)]=\frac{\mbox{d}}{\mbox{d}u}e(u), \nonumber \\
& & [d,f(u)]=\frac{\mbox{d}}{\mbox{d}u}f(u), \nonumber \\
& & [d,k^\pm_i(u)]=\frac{\mbox{d}}{\mbox{d}u}k^\pm_i(u), ~~i=1,2,
 \nonumber\\
& & k^\pm_i(u)k^\pm_j(v)=k^\pm_j(v)k^\pm_i(u),~~i,j=1,2, \nonumber\\
& & \rho(u_- -v_+)k^+_i(u)k^-_i(v)=k^-_i(v)k^+_i(u) \rho(u_+ -v_-),~~i=1,2,
 \nonumber\\
& & \rho(u_+ -v_- -\hbar)k^+_2(u)k^-_1(v)=k^-_1(v)k^+_2(u)
\rho(u_- -v_+ -\hbar), \nonumber\\
& & \rho(u_+ -v_- + \hbar)k^+_1(u)k^-_2(v)=k^-_2(v)k^+_1(u)
\rho(u_- -v_+ +\hbar), \nonumber\\
& & e(u)e(v)={u-v+\hbar\over u-v-\hbar}\ e(v)e(u),  \label{2.1}\\
& & f(u)f(v)={u-v-\hbar\over u-v+\hbar}\ f(v)f(u),  \nonumber\\
& & k^\pm_1(u)e(v)={u_\pm-v\over u_\pm-v+\hbar} e(v)k^\pm_1(u), \nonumber \\
& & k^\pm_2(u)e(v)={u_\pm-v\over u_\pm-v-\hbar} e(v)k^\pm_2(u), \nonumber \\
& & k^\pm_1(u)f(v)={u_\mp-v+\hbar\over u_\mp-v} f(v)k^\pm_1(u), \nonumber \\
& & k^\pm_2(u)f(v)={u_\mp-v-\hbar\over u_\mp-v} f(v)k^\pm_2(u), \nonumber \\
& & {[}e(u),f(v){]}=
{1\over h}\left(\delta(u_- -v_+)k^+_2(u_-)k^+_1(u_-)^{-1}
-\delta(u_+ -v_-)k^-_2(v_-)k^-_1(v_-)^{-1} \right), \nonumber
\end{eqnarray}

\noindent where

\begin{displaymath}
\delta(u-v)=\sum_{n+m=-1}u^nv^m, \hspace{1cm}
\delta(u-v)g(u)=\delta(u-v)g(v),
\end{displaymath}
\begin{displaymath}
u_\pm=u\pm{1 \over 4}\hbar c,
\end{displaymath}

\noindent and the function $\rho(u)$ is to be specified in the due course.
The equivalence of the Drinfeld currents to the Chevalley generators is
manifest in the following Laurent mode expansions of the currents,

\begin{eqnarray*}
e^{\pm}(u)=\pm\sum_{l\geq0\atop l<0} e[l] u^{-l-1},
f^{\pm}(u)=\pm\sum_{l\geq0\atop l<0} f[l] u^{-l-1},
k^{\pm}_i(u)=1\pm\hbar\sum_{l\geq0\atop l<0}k_i[l] u^{-l-1}.
\end{eqnarray*}

\noindent In the main context of this paper, we shall actually need the
equivalence between the Drinfeld currents and
Reshetikhin-Semenov-Tian-Shansky realization, the latter is given by the
following Yang-Baxter type relations \cite{RS},

\begin{eqnarray}
R^{\pm}(u-v) L^\pm_1(u) L^\pm_2(v) &=& L^\pm_2(v) L^\pm_1(u)
R^{\pm}(u-v),   \nonumber\\
R^+(u_- - v_+) L^+_1(u) L^-_2(v) &=& L^-_2(v) L^+_1(u) R^+(u_+ -v _-),
\label{2.2}
\end{eqnarray}

\noindent where

\begin{displaymath}
R^{\pm}(u)= \rho^{\pm}(u) \left(
\begin{array}{cccc}
1& & & \cr
 & $$u\over u+\hbar$$ & $$\hbar\over u+ \hbar$$ & \cr
 & $$\hbar\over u+\hbar$$ & $$u\over u+ \hbar$$ & \cr
 & & &1
\end{array}
\right),
\end{displaymath}
\begin{displaymath}
\rho^{\pm}(u)=\left({\Gamma^2\left({1\over 2} \mp {u\over 2\hbar}\right)}
\over{\Gamma\left(\mp {u\over 2\hbar}\right)
\Gamma\left(1\mp {u\over 2\hbar}\right)} \right)^{\pm 1},
\end{displaymath}

\noindent and the function $\rho(u)$ appeared in (\ref{2.1}) is precisely
$\rho^+(u)$. Notice that the scalar functions $\rho^{\pm}(u)$ in the
$R$ matrices $R^{\pm}(u)$ are chosen such that the unitarity and crossing
symmetry for the $R$ matrices hold, i.e.,

\begin{displaymath}
R^{+}(u)R^{-}(-u)=1,~~~~(C \otimes id)R^{\pm}(u)(C\otimes id)^{-1}
=R^{\mp}(-u-\hbar)^{t_1},
\end{displaymath}

\noindent where $t_1$ refers to the transpose in the first component space,
$C$ is the charge conjugation given in matrix form in the following,

\begin{displaymath}
C=\left(
\begin{array}{cc}
 & $$-1$$\cr
1&       \cr
\end{array}
\right).
\end{displaymath}

The equivalence between the two realizations (\ref{2.1}) and (\ref{2.2})
is an analog of the well known Ding-Frenkel equivalence \cite{DF}
of two similar realizations of $q$-affine algebras, see
also \cite{KL} in the case of Yangian
double $DY_\hbar(sl_2)_c$. In our case, the key point is that,
in eq.(\ref{2.2}), the $L^{\pm}(u)$ can be given a Gauss decomposition,

\begin{equation}
L^{\pm}(u)=
\left(\begin{array}{cc}
1&0\cr
$$\hbar f^{\pm}(u_\mp)$$&1
\end{array}\right)
\left(\begin{array}{cc}
$$k^\pm_1(u)$$&0\cr
0& $$k_2^{\pm}(u)$$
\end{array}\right)
\left(\begin{array}{cc}
1&$$\hbar e^{\pm}(u_\pm)$$\cr
0&1
\end{array}\right),  \label{2.a}
\end{equation}

\noindent where the diagonal entries $k^{\pm}_i(u)$ are identified with the
Drinfeld currents $k^\pm_i(u)$ in (\ref{2.1}), and the
off-diagonal entries $e^\pm(u)$
and $f^{\pm}(u)$ are related to the Drinfeld currents $e(u)$ and $f(u)$
by

\begin{displaymath}
e(u)=e^+(u)-e^-(u),\hspace{1cm} f(u)=f^+(u)-f^-(u).
\end{displaymath}

The algebra $DY_\hbar(gl_2)_c$ can be splitted into two subalgebras:
the Yangian double $DY_\hbar(sl_2)_c$ and a Heisenberg
subalgebra. The Heisenberg subalgebra is generated by the currents

\begin{equation}
K^{\pm}(u) \equiv k^\pm_2(u+\hbar) k^\pm_1(u) -1. \label{2.3}
\end{equation}

\noindent It is an easy practice to show that $K^{\pm}(u)$ actually commute
with all generating functions of $DY_\hbar(gl_2)_c$ and thus
generate a central subalgebra.
The Yangian double $DY_\hbar(sl_2)_c$ is thus obtained from
$DY_\hbar(gl_2)_c$ by taking the quotient with respect to this center.
The resulting generating relations differ from that of $DY_\hbar(gl_2)_c$
only in those involving the currents $k^\pm_i(u)$,

\begin{eqnarray*}
& & [d,h^\pm(u)]=\frac{\mbox{d}}{\mbox{d}u}h^\pm(u), \\
& & [h^\pm(u),~h^\pm(v)] = 0,\\
& & h^\pm(u)e(v)={u_\pm-v+\hbar\over u_\pm-v-\hbar} e(v)h^\pm(u), \\
& & h^\pm(u)f(v)={u_\mp-v-\hbar\over u_\mp-v+\hbar} f(v)h^\pm(u), \\
& & h^+(u)h^-(v)
={u_+ -v_- +\hbar\over u_- -v_+ +\hbar}
\cdot{u_- -v_+ -\hbar\over u_+ -v_- -\hbar} h^-(v)h^+(u),\\
\end{eqnarray*}

\noindent where $h^\pm(u)$ is defined as

\begin{equation}
h^\pm(u)=k^\pm_2(u)k^\pm_1(u)^{-1}. \label{2.4}
\end{equation}

In ending this section let us remark that one can recover the original
currents $k^\pm_{1,2}(u)$ from eq.(\ref{2.3}) and (\ref{2.4})
in the following form,

\begin{eqnarray}
& &k^+_1(u) = \prod_{l\geq 0} {h^+(u-(2l+1)\hbar) \over h^+(u-2l\hbar)},
\hspace{1cm}
k^+_2(u) = \prod_{l\geq 0} {h^+(u-(2l+1)\hbar) \over h^+(u-(2l+2)\hbar)},
 \nonumber\\
& &k^-_1(u) = \prod_{l\geq 0} {h^-(u+(2l+2)\hbar) \over h^-(u+(2l+1)\hbar)},
\hspace{1cm}
k^-_2(u) = \prod_{l\geq 0} {h^-(u+2l\hbar) \over h^-(u+(2l+1)\hbar)}.
\label{2.5}
\end{eqnarray}

\noindent These formulas will be used in Section 3.

\section{$\hbar$-deformed Sugawara construction}

Let

\begin{equation}
L(u) = L^-(u-{\hbar\over 2}) L^+(u+{\hbar\over 2})^{-1}.  \label{3.1}
\end{equation}

\noindent The trace

\begin{displaymath}
l(u) = tr L(u) = L_{11}(u) + L_{22}(u),
\end{displaymath}

\noindent as formal power series, would then lie in
the formal completion of $DY_\hbar(gl_2)_c$. Following \cite{RS}
we may conclude that at $c=-2$ the coefficients of $l(u)$ are central
elements of the formal completion of $DY_\hbar(gl_2)_c$.

To express $l(u)$ in terms of the equivalence between the
Drinfeld currents and the Reshetikhin-Semenov-Tian-Shansky
formalism, we will now set $c=-2$. It follows from (\ref{3.1}) and
(\ref{2.a}) that

\begin{eqnarray}
& &L_{11}(u) = k^-_1 (u-{\hbar\over 2} )
k^+_1 (u+{\hbar\over 2} )^{-1}
+\hbar^2 k^-_1 (u-{\hbar\over 2} ) e^+ (u )
k^+_2 (u+{\hbar\over 2} )^{-1}
f^+ (u+\hbar ) \nonumber\\
& &~~~~~~- \hbar^2 k^-_1 (u-{\hbar\over 2} ) e^-(u)
k^+_2 (u+{\hbar\over 2} )^{-1} f^+ (u+ \hbar ),  \nonumber\\
& &L_{22}(u) = k^-_2 (u-{\hbar\over 2} )
k^+_2 (u+{\hbar\over 2} )^{-1}
-\hbar^2 f^- (u-\hbar ) k^-_1 (u-{\hbar\over 2} ) e^+(u)
k^+_2 (u+{\hbar\over 2} )^{-1} \nonumber\\
& &~~~~~~+\hbar^2 f^- (u-\hbar )
k^-_1 (u-{\hbar\over 2} ) e^-(u)
k^+_2 (u+{\hbar\over 2} )^{-1}. \label{3.b}
\end{eqnarray}

The last two equations have not yet been written purely in
terms of Drinfeld currents. In order to do so, we have to combine
$e^\pm(u)$ into $e(u)$ and $f^\pm(u)$ into $f(u)$. The first step
will of cause be moving the
$f^+(u+\hbar)$ from the right of
$k^+_2 (u+{\hbar\over 2} )^{-1}$ to the left in $L_{11}(u)$,
$f^-(u-\hbar)$ from the left of
$k^-_1 (u-{\hbar\over 2} )^{-1}$ to the right in $L_{22}(u)$.
To achieve this we have to use the commutation relations for
$f^\pm(u)$ and $k^\pm_i(v)$. The required relations read

\begin{eqnarray*}
& & f^-(v_+) k^-_1(u) = {u-v\over u-v+\hbar} k^-_1(u) f^-(v_+)
+ {\hbar\over u-v+\hbar} f^-(u_+) k^-_1(u),\\
& & k^+_2(v)^{-1} f^+(u_-) ={u-v \over u-v+\hbar}
f^+(u_-) k^+_2(v)^{-1} + {\hbar\over u-v+\hbar}
k^+_2(v)^{-1} f^+(v_-).
\end{eqnarray*}

\noindent Multiplying by $u-v+\hbar$, we can see that, at points
$v=u+\hbar$, the last equations become

\begin{eqnarray*}
& & f^-(u_+) k^-_1(u)=k^-_1(u) f^-(u_+ + \hbar),\\
& & k^+_2(u)^{-1} f^+(u_-) = f^+(u_- - \hbar) k^+_2(u)^{-1}.
\end{eqnarray*}

\noindent Notice that, when $c=-2$, we have $u_\pm=u\mp \frac{1}{2}\hbar$.
Therefore, the above equations can be written as

\begin{eqnarray}
& & f^-(u-\hbar) k^-_1(u-\frac{1}{2}\hbar) =
 k^-_1(u-\frac{1}{2}\hbar) f^-(u), \nonumber\\
& & k^+_2(u+\frac{1}{2}\hbar)^{-1} f^+(u+\hbar) =
 f^+(u) k^+_2(u+\frac{1}{2}\hbar)^{-1}. \label{3.a}
\end{eqnarray}

\noindent Substituting (\ref{3.a}) into (\ref{3.b}), we are led to

\begin{eqnarray*}
& & L_{11}(u) = k^-_1(u-{\hbar\over 2}) k^+_1(u+{\hbar\over 2})^{-1}
+\hbar^2 k^-_1(u-{\hbar\over 2}) \left[ e^+(u) -
e^-(u) \right] f^+(u) k^+_2(u+{\hbar\over 2})^{-1}, \\
& & L_{22}(u) = k^-_2(u-{\hbar\over 2}) k^+_2(u+{\hbar\over 2})^{-1}
-\hbar^2 k^-_1(u-{\hbar\over 2}) f^-(u) \left[ e^+(u) -
e^-(u) \right] k^+_2(u+{\hbar\over 2})^{-1}.
\end{eqnarray*}

\noindent Finally, we have the expression for $l(u)$,

\begin{equation}
l(u)= k^-_1(u-{\hbar\over 2}) k^+_1(u+{\hbar\over 2})^{-1}
+ k^-_2(u-{\hbar\over 2}) k^+_2(u+{\hbar\over 2})^{-1}
+ \hbar^2 k^-_1(u-{\hbar\over 2}) ~:e(u) f(u):~
k^+_2(u+{\hbar\over 2})^{-1}, \label{3.c}
\end{equation}

\noindent where

\begin{equation}
:e(u)f(u): = e(u) f^+(u) - f^-(u) e(u).  \label{3.d}
\end{equation}

\noindent The final equation (\ref{3.c}) is just the
required $\hbar$-deformed Sugawara operator.

\section{Wakimoto module of $DY_\hbar(sl_2)_c$ and $\hbar$-deformed
Miura map}

Given the $\hbar$-deformed Sugawara construction, our next goal is
to show its connection to the corresponding (deformed) Miura map.
This can be fulfilled by making use of the Wakimoto module of
the Yangian double.

We shall adopt the Wakimoto module of $DY_\hbar(sl_2)_c$ given by
Konno \cite{KNO}\footnote{Our notations here defers from that of Konno
in \cite{KNO} in the following way: the boson $\lambda$ corresponds to
$a_\Phi$ of Konno, and $b$ and $c$ correspond to $a_\phi$ and $a_\chi$
respectively, and there is a shift of spectral parameters in the
bosonization formulas of Drinfeld current because our starting
definition of the Yangian double $DY_\hbar(sl_2)_c$ defers from that of
Konno by such a shift.}.

Introduce the following three sets of Heisenberg algebras with generators
$\lambda_n,~b_n,~c_n$ $n \in Z-\{0\}$, $\mbox{exp}(\pm q_\lambda),~
\mbox{exp}(\pm q_b),~\mbox{exp}(\pm q_c),~p_\lambda,p_b,$ and $p_c$,

\begin{eqnarray}
& &[\lambda_m,~\lambda_n]={k+2 \over 2} m \delta_{m+n,0},~~~~~~
[p_\lambda,~q_\lambda] ={k+2 \over 2}, \nonumber\\
& &[b_m,~b_n]=-m \delta_{m+n,0},~~~~~~~~~~[p_b,~q_b]=-1, \nonumber\\
& &[c_m,~c_n]= m \delta_{m+n,0},~~~~~~~~~~~~[p_c,~q_c]= 1. \label{4.1}
\end{eqnarray}

\noindent For $X=\lambda,b,c$, define

\begin{eqnarray*}
X(u;A,B)= \sum_{n>0} {X_{-n} \over n}(u+A\hbar)^n -
\sum_{n>0} {X_{n} \over n}(u+B\hbar)^{-n} + \mbox{log}(u+B\hbar) p_X +q_X,
\end{eqnarray*}

\noindent and together, $X(u;A)=X(u;A,A)$. We also use the abbreviations

\begin{eqnarray*}
& &X^+(u;B)= - \sum_{n>0} {X_{n} \over n}(u+B\hbar)^{-n} ,\\
& &X^-(u;A)= \sum_{n>0} {X_{-n} \over n}(u+A\hbar)^n .
\end{eqnarray*}

The triple-mode Fock space is defined as follows. Let $|0\rangle$
be a vector satisfying

\begin{displaymath}
X_n |0 \rangle =0,~~n>0;~~~~p_X | 0 \rangle =0.
\end{displaymath}

\noindent Then $|l,s,t\rangle \equiv \mbox{exp}
\left(\frac{l}{k+2}q_\lambda + s q_b + t q_c \right) | 0 \rangle$
is a vacuum state with $\lambda,b,c$ charges $l,-s,t$ respectively.
The Fock space is generated by the action of $\lambda_{-n},b_{-n},
~c_{-n} (n>0)$ on $|l,s,t\rangle$,

\begin{displaymath}
{\cal F}_{l,s,t}=\left\{ \prod_{n>0} \lambda_{-n}
\prod_{n'>0} b_{-n'} \prod_{n''>0} \right\} c_{-n''} | l,s,t \rangle.
\end{displaymath}

\noindent On ${\cal F}_{l,s,t}$, the normal ordering of
$\mbox{exp}(X(u;A,B))$ is defined as

\begin{eqnarray*}
:\mbox{exp}(X(u;A,B)):=\mbox{exp}(X^-(u;A)) \mbox{exp}(q_X)
(u+B\hbar)^{p_X}\mbox{exp}(X^+(u;B)).
\end{eqnarray*}

Now we are ready to define the Wakimoto module for the Yangian double
$DY_\hbar(sl_2)_c$. This is nothing but a homomorphism from the
above defined Heisenberg algebras to $DY_\hbar(sl_2)_c$ under which
the action of the Drinfeld currents on ${\cal F}_{l,s,t}$ is
given by \cite{KNO}

\begin{eqnarray*}
& & c=k,\\
& & d=d_\lambda+d_b+d_c,\\
& &~~~~~~ d_\lambda=\frac{2}{k+2}\left(\lambda_{-1}p_\lambda +
\sum_{n>0}\lambda_{-(n+1)}\lambda_n \right),\\
& &~~~~~~ d_b=-b_{-1}p_b- \sum_{n>0} b_{-(n+1)}b_n,\\
& &~~~~~~ d_c= c_{-1}p_c + \sum_{n>0} c_{-(n+1)} c_n,\\
& & h^+(u) = \mbox{exp}\left[\lambda^+(u;-\frac{3}{4}k)
- \lambda^+(u;-(\frac{3}{4}k+2)) \right.\\
& &~~~~~~\left. + b^+(u;-\frac{3}{4}k)-
b^+(u;-(\frac{3}{4}k+2))\right]
\left(\frac{u-\frac{3}{4}k\hbar}{u-(\frac{3}{4}k+2)\hbar}
\right)^{p_\lambda +p_b},\\
& & h^-(u) = \mbox{exp}\left[\frac{2}{k+2}
\left( \lambda^-(u;-(\frac{5}{4}k+3))
- \lambda^-(u;-(\frac{1}{4}k+1)) \right)\right] \\
& &~~~~~~ \times \mbox{exp}\left[ b^-(u;-(\frac{5}{4}k+3))
- b^-(u;-(\frac{1}{4}k+1)) \right],\\
& & e(u) = - \frac{1}{\hbar}: \left( \mbox{exp}(-c(u;-(k+1)))
-\mbox{exp}(-c(u;-(k+2))) \right) \\
& &~~~~~~ \times \mbox{exp}(-b(u;-(k+1),-(k+2))):, \\
& & f(u) = \frac{1}{\hbar} : \left(
\mbox{exp}\left[ \lambda^+(u;-\frac{1}{2}k) -
\lambda^+(u;-(\frac{1}{2}k+2)) \right]
\left(\frac{u-\frac{1}{2}k\hbar}{u-(\frac{1}{2}k+2)\hbar}
\right)^{p_\lambda} \right.\\
& &~~~~~~ \times \mbox{exp}\left[ b(u;-(\frac{1}{2}k+1),-\frac{1}{2}k)
+ c(u;-(\frac{1}{2}k+1)) \right] \\
& &~~~~~~ - \mbox{exp}\left[\frac{2}{k+2}
\left( \lambda^-(u;-(\frac{3}{2}k+3)) -
\lambda^-(u;-(\frac{1}{2}k+1)) \right) \right] \\
& &~~~~~~ \times \mbox{exp}\left.\left[ b(u;-(\frac{3}{2}k+3),
~-(\frac{3}{2}k+2)) + c(u;-(\frac{3}{2}k+2)) \right] \right):.
\end{eqnarray*}

Introducing the notations

\begin{eqnarray*}
& &\Upsilon^+(u)=(u-(\frac{1}{2}k+1)\hbar)^{p_\lambda}
\mbox{exp}\left(\lambda^+(u;-(\frac{1}{2}k+1))\right),\\
& &\Upsilon^-(u)=\mbox{exp}\left(\frac{2}{k+2}
\lambda^-(u;-(k+2)) \right),
\end{eqnarray*}

\noindent the expressions for $h^\pm(u)$ and $f(u)$ can be recasted into
a relatively shorter form,

\begin{eqnarray}
& & h^+(u) = \Upsilon^+(u_- + \hbar) \Upsilon^+(u_- -\hbar)^{-1}
\left(\frac{u-\frac{3}{4}k\hbar}{u-(\frac{3}{4}k+2)\hbar}
\right)^{p_b}  \nonumber\\
& &~~~~~~\times \mbox{exp}\left[b^+(u;-\frac{3}{4}k)-
b^+(u;-(\frac{3}{4}k+2))\right], \label{a.a} \\
& & h^-(u) = \Upsilon^-(u_+ -\frac{k+2}{2}\hbar)
\Upsilon^-(u_+ + \frac{k+2}{2}\hbar)^{-1}  \nonumber\\
& &~~~~~~\times \mbox{exp}\left[ b^-(u;-(\frac{5}{4}k+3))
- b^-(u;-(\frac{1}{4}k+1)) \right], \label{a.b} \\
& & f(u) = \frac{1}{\hbar}
: \Upsilon^+(u + \hbar) \Upsilon^+(u -\hbar)^{-1}  \nonumber \\
& &~~~~~~\times \mbox{exp}\left[ b(u;-(\frac{1}{2}k+1),-\frac{1}{2}k)
+ c(u;-(\frac{1}{2}k+1)) \right]  \nonumber\\
& &~~~~~~ - \Upsilon^-(u -\frac{k+2}{2}\hbar)
\Upsilon^-(u + \frac{k+2}{2}\hbar)^{-1} \nonumber \\
& &~~~~~~\times \mbox{exp}\left[ b(u;-(\frac{3}{2}k+3),
~-(\frac{3}{2}k+2)) + c(u;-(\frac{3}{2}k+2)) \right]:.  \nonumber
\end{eqnarray}

In order to obtain the $\hbar$-deformed Miura map, we needs to express
$l(u)$ in terms of Laurent modes of only one of the three bosons, $\lambda$.
To achieve this goal, we need an bosonic expression for $k^\pm_i(u)$, which
can be obtained by substituting (\ref{a.a}) and (\ref{a.b}) into (\ref{2.5}),

\begin{eqnarray}
& &k^+_1(u) = \frac{\Upsilon^+\left(u-\frac{k}{2}\hbar \right)}
{\Upsilon^+\left(u-\frac{k}{2}\hbar+\hbar\right)}
\left(\frac{u-(k+1)\hbar}{u-k\hbar}\right)^{p_b} \nonumber\\
& &~~~~~~~\times \mbox{exp}\left(b^+(u;-(k+1))-b^+(u;-k)\right), \nonumber\\
& &k^-_1(u)= \prod_{l \geq 0}
\frac{\Upsilon^-\left(u+k\hbar+(2l+2)\hbar\right)
\Upsilon^-\left(u+(2l+1)\hbar\right)}
{\Upsilon^-\left(u+k\hbar+(2l+3)\hbar\right)
\Upsilon^-\left(u+2l)\hbar\right)} \nonumber\\
& &~~~~~~~\times \mbox{exp}\left(b^-(u;-(k+1))-b^-(u;-(k+2)\right),
\nonumber\\
& &k^+_2(u) = \frac{\Upsilon^+\left(u-\frac{k}{2}\hbar \right)}
{\Upsilon^+\left(u-\frac{k}{2}\hbar-\hbar\right)}
\left(\frac{u-(k+1)\hbar}{u-(k+2)\hbar}\right)^{p_b} \nonumber\\
& &~~~~~~~\times \mbox{exp}\left(b^+(u;-(k+1))-b^+(u;-(k+2)\right)  ,
\nonumber\\
& &k^-_2(u)= \prod_{l \geq 0}
\frac{\Upsilon^-\left(u+k\hbar+(2l+2)\hbar\right)
\Upsilon^-\left(u+(2l-1)\hbar\right)}
{\Upsilon^-\left(u+k\hbar+(2l+1)\hbar\right)
\Upsilon^-\left(u+2l)\hbar\right)}  \nonumber\\
& &~~~~~~~\times \mbox{exp}\left(b^-(u;-(k+3))-b^-(u;-(k+2)\right).
\label{b.1}
\end{eqnarray}

On the other hand, on the Fock space ${\cal F}_{l,s,t}$, the normal
ordering for $e(u)$ and $f(u)$ is given by

\begin{equation}
:e(u)f(u): = - \int_{C_1} \mbox{d}v\frac{e(u)f(v)}{u-v-\frac{k+2}{2}\hbar}
+ \int_{C_2} \mbox{d}v\frac{f(v)e(u)}{u-v+\frac{k+2}{2}\hbar},
\label{a.c}
\end{equation}

\noindent where $C_1$ and $C_2$ are respectively circles of radius
$|v| > |u-\frac{k+2}{2}\hbar|$ and $|v| < |u+\frac{k+2}{2}\hbar|$.
Applying the concrete expressions for $e(u)$ and $f(v)$ in
(\ref{a.c}), we can get

\begin{eqnarray}
& & :e(u)f(u): = \frac{1}{\hbar^2} (
- \frac{\Upsilon^+\left(u-(\frac{k}{2}-1)\hbar \right)}
{\Upsilon^+\left(u-(\frac{k}{2}+1)\hbar\right)}  \nonumber \\
& &~~~~~~~~\times:\mbox{exp}\left(b(u;-(k+1),-k)-b(u;-(k+1),-(k+2)\right):
\nonumber\\
& &~~~~~~+ \frac{\Upsilon^+\left(u-(\frac{k+2}{2}-1)\hbar \right)}
{\Upsilon^+\left(u-(\frac{k+2}{2}+1)\hbar\right)}  \nonumber\\
& &~~~~~~~~\times :\mbox{exp}\left(b(u;-(k+2),-(k+1))
-b(u;-(k+1),-(k+2)\right):  \nonumber\\
& &~~~~~~+ \frac{\Upsilon^-\left(u\right)}
{\Upsilon^-\left(u+(k+2)\hbar\right)}  \nonumber\\
& &~~~~~~~~\times :\mbox{exp}\left(b(u;-(k+2),-(k+1)
-b(u;-(k+1),-(k+2)\right): \nonumber \\
& &~~~~~~- \frac{\Upsilon^-\left(u-\hbar \right)}
{\Upsilon^+\left(u-(k+1)\hbar\right)}  \nonumber\\
& &~~~~~~~~\times :\mbox{exp}(b(u;-(k+3),-(k+2)-b(u;-(k+1),-(k+2) ): )
\label{b.2}
\end{eqnarray}

\noindent Substituting equations (\ref{b.1})-(\ref{b.2}) into
(\ref{3.c}) and after some
some algebra, we finally obtain

\begin{equation}
l(u)= \Lambda(u-\frac{\hbar}{2}) + \Lambda(u+\frac{\hbar}{2})^{-1},
\label{Mu}
\end{equation}

\noindent where

\begin{equation}
\Lambda(u) = \Lambda^+(u) \Lambda^-(u) \label{La}
\end{equation}

\noindent and

\begin{eqnarray*}
& & \Lambda^+(u) =
\frac{\Upsilon^+\left(u-\frac{k+2}{2}\hbar + \frac{\hbar}{2}\right)}
{\Upsilon^+\left(u-\frac{k+2}{2}\hbar - \frac{\hbar}{2}\right)}, \\
& & \Lambda^-(u) =
\prod_{l \geq 0}
\frac{\Upsilon^-\left(u+(k+2)\hbar+(2l+1)\hbar - \frac{\hbar}{2}\right)
\Upsilon^-\left(u+(2l+1)\hbar + \frac{\hbar}{2}\right)}
{\Upsilon^-\left(u+(k+2)\hbar+(2l+1)\hbar+ \frac{\hbar}{2} \right)
\Upsilon^-\left(u+(2l+1)\hbar - \frac{\hbar}{2} \right)}.
\end{eqnarray*}

\noindent Notice that eq.(\ref{Mu}) contains only expressions involving
the Laurent modes of $\lambda$ and thus provides a map from one
free $\hbar$-deformed bosonic field $\lambda$ to the $\hbar$-deformed
Sugawara operator $l(u)$. This is precisely the desired $\hbar$-deformed
Miura map.

We have to remark that, though in the form of eq.(\ref{Mu})
the $\hbar$-deformed
Miura map looks very similar to the $q$-deformed version, the
actual way of mapping from the bosonic field $\lambda$ to $l(u)$
is much more complicated than the $q$-deformed case. The
complexity comes about in the infinite product structure in the expression
for $\Lambda^-(u)$. However, despite such complexities the
operator product between $\Lambda^+$ and $\Lambda^-$ is rather simple.
It reads

\begin{displaymath}
\Lambda^+(u) \Lambda^-(v) = \frac{\rho(u-v-(k+2)\hbar)}{\rho(u-v)}
\Lambda^-(v) \Lambda^+(u).
\end{displaymath}

Another remark is concerned with the fact that the Sugawara operator $l(u)$ is
associated with the two-dimensional representation of $sl_2$.
To obtain analogous operators associated with higher dimensional
representations of $sl_2$, one can apply the fusing procedure which is
similar to the $q$-deformed case, i.e. the operator $l^{(n)}(u)$ associated
with the $(n+1)$-dimensional representation of $sl_2$
can be obtained from the following iterative relation,

\begin{eqnarray*}
& &l^{(1)}(u-n\hbar) l^{(n)}(u) = l^{(n+1)}(u) + l^{(n-1)}(u),\\
& &l^{(1)}(u) = l(u),~~~~~l^{(0)}(u)=1.
\end{eqnarray*}

\noindent The explicit form for $l^{(n)}(u)$ reads

\begin{eqnarray*}
& &l^{(n)}(u) =
\Lambda(u-\frac{\hbar}{2}) \Lambda(u-\frac{3 \hbar}{2})
\Lambda(u-\frac{5 \hbar}{2}) ... \Lambda(u-\frac{(2n-1) \hbar}{2})\\
& &~~~~~~~ +
\Lambda(u+\frac{\hbar}{2})^{-1} \Lambda(u-\frac{3 \hbar}{2})
\Lambda(u-\frac{5 \hbar}{2}) ... \Lambda(u-\frac{(2n-1) \hbar}{2})\\
& &~~~~~~~ +
\Lambda(u+\frac{\hbar}{2})^{-1} \Lambda(u-\frac{\hbar}{2})^{-1}
\Lambda(u-\frac{5 \hbar}{2}) ... \Lambda(u-\frac{(2n-1) \hbar}{2})\\
& &~~~~~~~ +
\Lambda(u+\frac{\hbar}{2})^{-1} \Lambda(u-\frac{\hbar}{2})^{-1}
\Lambda(u-\frac{3 \hbar}{2})^{-1} ... \Lambda(u-\frac{(2n-1) \hbar}{2})\\
& &~~~~~~~ + ......\\
& &~~~~~~~ +
\Lambda(u+\frac{\hbar}{2})^{-1} \Lambda(u-\frac{\hbar}{2})^{-1}
\Lambda(u-\frac{3 \hbar}{2})^{-1} ... \Lambda(u-\frac{(2n-3) \hbar}{2})^{-1}.
\end{eqnarray*}

\section{Poisson brackets for the $\hbar$-deformed Virasoro algebra}

Let us recall that in the undeformed case, the Miura map provides
a free field representation of the classical (i.e. Poisson bracket)
Virasoro algebra. In the $q$-deformed case, such a map also gives
rise from the Poisson brackets for a $q$-deformed bosonic field
to a $q$-deformed Virasoro Poisson algebra. In this section we shall
show how we can obtain an analogous $\hbar$-deformed algebra.

First let us explain how a Poisson brackets could arise from a
purely quantum theory. The key point is as follows. When
$k+2$ approaches zero, the commutation relations in (\ref{4.1}), divided by
$k+2$, naturally induces a Poisson bracket structure,

\begin{equation}
\{ \lambda_m,~\lambda_n \} = \frac{1}{2}m \delta_{m+n,0},~~~~~
\{ p_\lambda,~q_\lambda \} = \frac{1}{2}. \label{5.1}
\end{equation}

\noindent This Poisson structure turn the quantum $\hbar$-deformed
free field $\lambda(u;A,B)$ into a classical object. In the meantime,
all functions of the quantum field $\lambda$ are also turned into
classical ones, i.e. the noncommuting objects become now commutative
and the original commutation relations are turned into
Poisson brackets.

To obtain the Poisson bracket for the $\hbar$-deformed Virasoro algebra,
we need to take the limit of $l(u)$ as $k+2 \rightarrow 0$. Simply
substituting $k+2=0$ into the expressions of $\Upsilon^\pm(u)$ does not make
sense because $\Upsilon^-(u)$ is not a well defined object
as $k+2 \rightarrow 0$. However, at level of $\Lambda^\pm(u)$, well defined
limits could be obtained. The limits of $\Lambda^\pm(u)$ read

\begin{eqnarray*}
& & \Lambda^\pm(u) = A^\pm(u-\frac{1}{2} \hbar)
A^\pm(u+\frac{1}{2} \hbar)^{-1},\\
& & A^+(u) = u^{-p_\lambda}
\mbox{exp}\left\{ \sum_{l>0}
\frac{\lambda_n}{n} u^{-n} \right\},\\
& & A^-(u) = \prod_{l \geq 0} B^-(u+(2l+1)\hbar),\\
& & B^-(u) = \mbox{exp}\left\{\sum_{n>0} 2\hbar
\lambda_{-n} u^{n-1} \right\}.
\end{eqnarray*}

Using the Poisson brackets (\ref{5.1}), one can easily calculates

\begin{eqnarray*}
& & \{ A^+(u),~B^-(v) \} = \hbar \frac{\partial}{\partial v}
\sum_{n>0} \frac{1}{n} \left(\frac{v}{u} \right)^n
A^+(u) B^-(v)\\
& &~~~~~~~= - \hbar \frac{\partial}{\partial v}
\mbox{log} \left(1- \frac{v}{u} \right)
A^+(u) B^-(v), ~~~~~(|u| > |v|) \\
& & \{ B^-(u),~A^+(v) \} = \hbar \frac{\partial}{\partial u}
\mbox{log} \left(1- \frac{u}{v} \right)
B^-(u) A^+(v),  ~~~~~(|v| > |u|)
\end{eqnarray*}

\noindent from which we obtain

\begin{eqnarray*}
& & \{ A^+(u),~A^-(v) \} =
- \hbar \frac{\partial}{\partial v}
\sum_{l \geq 0} \mbox{log} \left(1- \frac{v+(2l+1)\hbar}
{u} \right) A^+(u) A^-(v) \\
& &~~~~~~~ = - \hbar \frac{\partial}{\partial v}
\mbox{log} \left[ \prod_{l \geq 0} \left(1- \frac{v+(2l+1)\hbar}
{u} \right) \right] A^+(u) A^-(v), \\
& & \{ A^-(u),~A^+(v) \} =
\hbar \frac{\partial}{\partial u}
\mbox{log} \left[ \prod_{l \geq 0} \left(1- \frac{u+(2l+1)\hbar}
{v} \right) \right] A^-(u) A^+(v),
\end{eqnarray*}

\noindent and further,

\begin{eqnarray*}
& &\{ \Lambda^+(u),~\Lambda^-(v) \} = \hbar \frac{\partial}{\partial v}
\mbox{log} \rho(u-v) \Lambda^+(u) \Lambda^-(v), \\
& &\{ \Lambda^-(u),~\Lambda^-(u) \} = - \hbar \frac{\partial}{\partial u}
\mbox{log} \rho(v-u) \Lambda^-(u) \Lambda^+(v).
\end{eqnarray*}

\noindent Remembering the definition (\ref{La}) of $\Lambda(u)$, we have

\begin{displaymath}
\{ \Lambda(u),~\Lambda(v) \} =
\hbar \left[ \frac{\partial}{\partial v}
\mbox{log} \rho(u-v) - \frac{\partial}{\partial u}
\mbox{log} \rho(v-u) \right] \Lambda (u) \Lambda(v).
\end{displaymath}

\noindent Finally, we have for

\begin{displaymath}
l(u) \rightarrow s(u) = \Lambda(u-\frac{1}{2}\hbar) +
\Lambda(u+\frac{1}{2}\hbar)^{-1}
\end{displaymath}

\noindent the following Poisson bracket,

\begin{eqnarray*}
& & \{ s(u),~s(v) \}  =
\hbar \left[ \frac{\partial}{\partial v}
\mbox{log} \rho(u-v) - \frac{\partial}{\partial u}
\mbox{log} \rho(v-u) \right] s(u) s(v) \\
& &~~~~~~~ + \hbar \delta(u-v-\hbar) - \hbar \delta(u-v+\hbar),
\end{eqnarray*}

\noindent where the $\delta$-functions are defined in the same way as
in eq.(\ref{2.1}).

Before concluding this paper let us consider the connections between
the $\hbar$-deformed Miura map and some $\hbar$-difference equations.
These are analogs of the well-known Gelfand-Dickey equations written in the
simplest case: the classical Miura map

\begin{equation}
\partial^2 - q(t) = \left( \partial - \frac{1}{2}\chi(t) \right)
\left( \partial + \frac{1}{2}\chi(t) \right). \label{xx}
\end{equation}

\noindent Let ${\cal D}_\hbar$ be the $\hbar$-``derivative'' defined as

\begin{eqnarray*}
{\cal D}_\hbar Q(u) = Q(u-\hbar).
\end{eqnarray*}

\noindent Then the $\hbar$-deformed version of eq.(\ref{xx})
can be written as

\begin{eqnarray*}
\left(\Lambda(u-\frac{\hbar}{2}) {\cal D}_\hbar -1 \right)
\left(\Lambda(u+\frac{\hbar}{2})^{-1} {\cal D}_\hbar -1 \right)
={\cal D}_\hbar^2 - s(u) {\cal D}_\hbar +1.
\end{eqnarray*}

\noindent In particular, if $Q(u)$ is a solution of the $\hbar$-difference
equation

\begin{eqnarray*}
\left({\cal D}_\hbar^2 - s(u) {\cal D}_\hbar +1 \right)Q(u+\hbar)=0,
\end{eqnarray*}

\noindent then $s(u)$ can be expressed in the form

\begin{eqnarray*}
s(u)= \frac{Q(u-\hbar)}{Q(u)} + \frac{Q(u+\hbar)}{Q(u)}.
\end{eqnarray*}

\noindent This equation is completely in analogy to the $q$-deformed version
given in Ref.\cite{FR}, which was first used by Baxter
\cite{baxter} in studying eight vertex
model. Such kind of equations have close connections to Bethe ansatz
equations.

\section{Concluding remarks and out-looking}

In this paper we constructed the $\hbar$-deformed Miura map and the
corresponding Virasoro algebra. This algebra can be viewed as either
a Yangian deformation of the classical Virasoro algebra or a
scaling limit of the $q$-deformed Virasoro algebra obtained by Frenkel
and Reshetikhin in \cite{FR}. However, the construction in this paper is much
more complicated than the $q$-deformed case because the $\hbar$-deformed
Sugawara operator involves an infinite product structure in terms of the
bosonic expression $\Upsilon^-(u)$. As a classical (Poisson bracket) algebra,
our algebra is the classical limit of the quantum $\hbar$-deformed
Virasoro algebra obtained from the quantum $q$-deformed Virasoro algebra
by taking proper limit, and may also be viewed as a deformation of the
conventional quantum Virasoro algebra governing the conformal field theory.

It is desirable that there exists an infinite family of $\hbar$-deformed
algebras---call them $\hbar$-deformed $W$-algebras---each corresponds
to a Yangian double of different underlying Lie algebras.
The $q$-deformed $W$-algebras are already constructed in \cite{FR} and
\cite{SKAO,LUKU,AKOS,AKOS2,FF}
both in classical and quantum form. However the $\hbar$-deformation
of $W$-algebras is not known yet and we hope to work out this problem in our
next publication.

Perhaps the most important application of the conventional classical
$W$-algebras (i.e. Gelfand-Dickey Poisson algebras) is in the Hamiltonian
description of integrable hierarchies such as the KdV hierarchy. For
$q$-deformed $W$-algebras the corresponding differential-difference
systems as deformed integrable hierarchies were obtained in \cite{frenkel}.
It is interesting to perform the analogous constructions for
$\hbar$-deformed $W$-algebras.

In the conventional (undeformed) cases, $W$-algebra generators are
connected to principal co-minors of certain Wronskian determinants.
It is quite interesting to ask the question that whether the $q$- and/or
$\hbar$-analogs exist. Such analogies are important if we want to identify
the existence of the $q$- and/or $\hbar$-versions of nonstandard
$W$-algebras, namely $W$-algebras beyond the standard $W_n$ series,
$W_3^{(2)}$ for instance. We also hope to consider these problems in
our future studies.

\section*{Acknowledgement}

One of the authors (Ding) would like to thank Professors
Ke Wu and Shi-Kun Wang for valuable discussions.

\end{document}